\documentclass[11pt,twoside]{article}

\usepackage[utf8]{inputenc}
\usepackage[T2A]{fontenc}
\usepackage[english]{babel}
\usepackage{amsmath, amssymb}
\usepackage{geometry}
\usepackage{setspace}
\usepackage{tikz}
\usepackage{fancyhdr}

\newcommand{\volume}{}
\newcommand{\issue}{}
\newcommand{\pubyear}{}

\renewcommand{\volume}{22}
\renewcommand{\issue}{1}
\renewcommand{\pubyear}{2025}

\newtheorem{lemma}{Lemma}
\newtheorem{definition}[lemma]{Definition}

\newtheorem{theorem}[lemma]{Theorem}

\def\endproof{\hfill$\Box$}

\fancypagestyle{plain}{
    \fancyhf{}

    \fancyhead[RE]{\small HERALD OF THE KAZAKH-BRITISH\\ TECHNICAL UNIVERSITY}
    \fancyhead[LE]{\small Vol.~\volume, No.~\issue, \pubyear}

    \fancyfoot[LO]{\thepage}
    \fancyfoot[RE]{\thepage}
    
}

\makeatletter

\renewcommand\@biblabel[1]{#1.}
\renewcommand{\thebibliography}[1]{%
  \section*{}%
  \begin{center}
    \textbf{\normalsize REFERENCES}
  \end{center}
  \list{\@biblabel{\arabic{enumiv}}}{%
    \settowidth\labelwidth{\@biblabel{#1}}%
    \leftmargin\labelwidth
    \advance\leftmargin\labelsep
    \usecounter{enumiv}%
  }%
  \sloppy
}
\makeatother

\begin{document}
\setcounter{page}{223}
\thispagestyle{plain}
\noindent
UDC 510.67 \\ 
IRSTI 27.03.66 \\ 

\vspace{1em}

\begin{center}
$^{1}$ Ussentay A.B., \\
PhD student, ORCID ID: 0009-0009-9036-2323, \\
e-mail: baglambekkyzy\_a@mail.ru

\vspace{0.5em}

\noindent
$^{2}$ *Sadykova K.K, \\
PhD, ORCID ID:  0009-0000-8651-312X, \\
*e-mail: sadkelbet@gmail.com
\vspace{0.5em}

\noindent
$^{1}$ L.N. Gumilyov Eurasian National University, Astana, Kazakhstan \\
$^{1}$ Geometry LLP , Astana, Kazakhstan \\
$^{2}$ S.Seifullin Kazakh Agro Technical Research University, Astana, Kazakhstan

\vspace{1em}

\textbf{O’Neil-type inequalities in Morrey spaces on the set $\mathbb{T}^n$}

\vspace{1em}

\textbf{Abstract}
\end{center}
\vspace{-0.5em}
\begin{small}

This research paper is devoted to the study of O’Neil-type inequalities for operator convolution in Morrey spaces. As Morrey-type spaces allow for a more precise characterisation of the local properties of functions, they play an important role in the theory of differential equations and harmonic analysis. This article analyses the development of the classical Minkowski and Young inequalities in the context of convolution, and examines their generalisation from Lorentz spaces to Morrey spaces.  The main aim of the study is to determine the conditions for the boundedness of O’Neill-type inequalities in local and generalised Morrey spaces defined on the $n$-dimensional unit cube. This paper presents a more detailed analysis of the classical Minkowski and Young inequalities. In addition to these inequalities, it generalises the well-known O’Neil inequality, defined for the Lorentz space. The paper focuses on proving interpolation embeddings and establishing relationships between various functional structures via dyadic decompositions. This work focuses on proving interpolation embeddings and establishing relationships between various functional structures through dyadic decompositions. The main result of this research is the proof that the density operators under consideration are bounded in functional spaces subject to certain parametric constraints.  New estimates have been obtained for smooth functions in the $n$-dimensional unit cube. Furthermore, interpolation embeddings and relationships between local and generalised Morrey spaces have been established.

\vspace{0.5em}

\end{small}

\vspace{1em}

\textbf{Introduction}

Let $\mathbb{T}^n=[0,1]^n$ be the $n$-dimensional unit cube. Let
$f$ and $g$ be functions defined on $\mathbb{T}^n$. The convolution of two functions  $f*g$ is defined by:
$$(f*g)(x)=\int_{\mathbb{T}^n}f(x-y)g(y)dy,$$
whenever the integral exists. This operation combines two functions to produce a new function that describes their averaging properties. The convolution operator plays a fundamental role in harmonic analysis, partial differential equations, and signal processing. 

One of the first fundamental inequalities for convolution is Minkowski's inequality \cite{Gl2014} :
Let $f\in L_p$ and $g\in L_1$, where  $1\le p\le \infty$. Then
$$\|f*g\|_{L_p}\le \|f\|_{L_p}\|g\|_{L_1}.$$

This result shows that convolution with a function from $L_1$ defines a bounded operator on $L_p$. A more general estimate is provided by Young's inequality \cite{Bennet&Sharley}. Let $1\le p,q,r\le \infty$ satisfy
$$\frac{1}{p}+\frac{1}{q}=1+\frac{1}{r}$$.
Then $f\in L_p$ and $g\in L_q$ convolution of the functions satisfies the following inequality:
$$\|f*g\|_{L_r}\le \|f\|_{L_p}\|g\|_{L_q}.$$
This inequality generalises many classical convolution estimates.

The Riesz potential is one of the fundamental operators in harmonic analysis and the theory of partial differential equations. For $\alpha< n$ it is defined by:
$$I_{\alpha}f(x)=\int_{\mathbb{T}^n}\frac{f(y)}{|x-y|^{n-\alpha}}dy.$$
If $1<\alpha<n$, $1<p<\frac{n}{\alpha}$, $0<\lambda<n$  satisfy and   $\frac{1}{q}=\frac{1}{p}-\frac{\alpha}{n}$, then the following estimate holds
$$\|I_{\alpha}f\|\le C\|f\|_{M_p^{\lambda}}.$$

A significant development in the theory of convolution inequalities appeared in the work of O'Neil \cite{Conv 1963}.who extended Young's inequality to Lorentz spaces  $L_{p,r}$ Lorentz spaces generalise Lebesgue spaces and provide a more refined description of integrability and function distribution. Consequently, they are particularly useful in the study of limiting cases and interpolation phenomena.
Let $f\in L_{p_1,r_1}$, $g\in L_{p_2,r_2}$ and $1\le p,p_1,p_2<\infty$, $0<r,r_1,r_2\le \infty$. Assume that
$$\frac{1}{p_1}+\frac{1}{p_2}=1+\frac{1}{p}, \,\,\frac{1}{r_1}+\frac{1}{r_2}=\frac{1}{r}$$. 
then $f*g\in L_{p,r}$ the following estimate holds
$$\|f*g\|_{L_{p,r}}\le \|f\|_{L_{p_1,r_1}}\|g\|_{L_{p_2,r_2}}.$$

 While Young's inequality provides boundedness in Lebesgue spaces $L_p$, O'Neil's result yields a more refined estimate in Lorentz spaces. In particular, it distinguishes functions having the same Lebesgue $L_p$ norm but different distribution properties. Therefore, O'Neil's inequality has become an important tool in interpolation theory and the study of singular integral operators

A significant advancement in the theory of convolution inequalities was achieved in \cite{Tleu Sadyk19},\cite{TleuSad2020}  where O'Neil's inequality was extended to anisotropic Lorentz spaces.

Building upon O'Neil's result, E. Nursultanov and S. Tikhonov determined the extremal and critical conditions for O'Neil-type embedding inequalities in Lorentz spaces. In \cite{NursTikh} these inequalities were investigated comprehensively. The authors established not only the exact range of validity of these inequalities but also identified the cases in which they fail to hold.

In \cite{ED&SD} E. Nursultanov and D. Suragan generalised O'Neil's inequality to Morrey-type spaces. Their work also produced new results concerning the boundedness of the Riesz potential operator in Morrey-type spaces. These results provide sharper estimates for singular operators and broaden the functional-analytic framework for studying local regularity properties.




\textbf{Materials and Methods}

\begin{definition}
Let $0 \le \lambda \le \frac{n}{p}$ and $0<p<\infty$. The set of all functions $f \in L_p^{\mathrm{loc}}(\mathbb{T}^n)$ is called the \emph{Morrey space} if the norm 
\[
\|f\|_{M_p^\lambda(\mathbb{T}^n)} 
= \sup_{x \in \mathbb{T}^n} \sup_{r>0} r^{-\lambda} \|f\|_{L_p(Q_r(x))} < \infty
\]
where $Q_r(x)$ — denotes the cube centered at $x$ with radius $r>0$.

In the case $\lambda=0$, we have
\[
M_p^0(\mathbb{T}^n) = L_p(\mathbb{T}^n),
\]
Indded, if $f \in M_p^0$, then
\[
\|f\|_{M_p^0} = \sup_{x \in \mathbb{T}^n} \sup_{r>0} r^{0} \|f\|_{L_p(Q_r(x))} = \|f\|_{L_p(\mathbb{T}^n)}.
\]

If $\lambda = \frac{n}{p}$ and $0<p<\infty$, then
\[
M_p^{\frac{n}{p}}(\mathbb{T}^n) = L_\infty(\mathbb{T}^n).
\]
\end{definition}
\begin{definition}
   Let $k=0,1,...,\infty$. Denote by $G_k$ the collection of all dyadic cubes in  $\mathbb{T}^n$ of the form:
\[
[0,2^{-k})^n + 2^{-k} m, \quad m=(m_1,m_2,...,m_i), \quad m_i=0,1,..,2^{(k-1)}, i=\overline{1,k}.
\]
Then it is clear that
\[
\mathbb{T}^n = \bigsqcup_{Q \in G_k} Q
\]
where $\bigsqcup$ — denotes the disjoint union of sets. 

Let the family of all dyadic cubes in $\mathbb{T}^n$ be defined by
\[
G = \bigcup_{k \in \mathbb{Z}} G_k.
\]

Each cube $Q \in G_k$  is partitioned into $2^n$  cubes belonging to $G_{k-1}$

A set of  disjoint cubes $\mathfrak{T} = \{Q\} \subset G$ is called a dyadic partition of $\mathbb{T}^n$ if:

    $$\mu (\mathbb{T}^n\backslash{\displaystyle\bigcup_{Q \in \mathfrak{T}} Q})=0.$$
\end{definition}
\begin{definition}
Let $k=0,1,...,\infty$, $\lambda \in \mathbb{R}$, $0<p,q \le \infty$ and let $\mathfrak{T}=(Q)$ be a dyadic partition of $\mathbb{T}^n$ .   We define the local Morrey space$LM_{p,q}^\lambda(\mathfrak{T})$ as the set of all measurable functions $f$ such that
\[
\|f\|_{LM_{p,q}^\lambda{(\mathfrak{T})}} 
= \Biggl( \sum_{k \in \mathbb{N}} \Bigl( 2^{-k \lambda} \sum_{Q \in \mathfrak{T}_k = \mathfrak{T} \cap G_k} \|f\|_{L_p(Q)} \Bigr)^q \Biggr)^{1/q}<\infty.
\]
 
\end{definition}
\begin{definition}
Let $\Omega \subset \mathbb{T}^n$, $0<p<\infty$, $0<q\le \infty$, $0\le \lambda \le \frac{n}{p}$ . We define the generalized Morrey-type space$M_{p,q}^\lambda(\Omega)$  as the collection of all Lebesgue measurable functions $f \in L_p^{\mathrm{loc}}(\mathbb{T}^n)$ such that
\[
\|f\|_{M_{p,q}^\lambda(\Omega)} 
= \left( \sum_{k=0}^\infty \bigl( 2^{-k\lambda} \sup_{Q \in G_k} \|f\|_{L_p(Q)} \bigr)^q  \right)^{1/q}<\infty.
\]
he spaces defined above are generalizations of the classical Morrey spaces. Indeed, when $q=\infty$, one recovers the classical Morrey space:   
\[
M_{p,\infty}^\lambda = M_p^\lambda.
\]  

\end{definition}

\vspace{1em}

\textbf{Results and Discussion}

\begin{lemma} \cite{Bur.Chigam}  
\begin{enumerate}
    \item[(i)] If $1 \le p_0 < p_1 < \infty$, then
    \[
    M_{p_1,q}^{\alpha_1} \hookrightarrow M_{p_0,q}^{\alpha_0}, \quad \text{where } \alpha_0 = \alpha_1 - \frac{n(p_0 - p_1)}{p_0 p_1}.
    \]
    \item[(ii)] If $0 < q_0 < q_1 \le \infty$, then
    \[
    M_{p,q_0}^\alpha \hookrightarrow M_{p,q_1}^\alpha.
    \]
\end{enumerate}  
\end{lemma}

\begin{theorem}\label{th7} \cite{ED&SD} 
Let $\lambda_0 \neq \lambda_1$, $0<p\le \infty$, $0<q,q_0,q_1\le \infty$  $\theta \in (0,1)$, then next inequality holds
\[
(M_{p,q_0}^{\lambda_0}, M_{p,q_1}^{\lambda_1})_{\theta,q} \hookrightarrow M_{p,q}^\lambda,
\]
where
\[
\lambda = (1-\theta)\lambda_0 + \theta \lambda_1.
\]

\end{theorem}

\emph{The Köthe dual space} (or simply the dual space)   of the generalized Morrey space $M_{p,q}^\lambda$   is denoted by $ (M_{p,q}^\lambda)'$ and is defined by
\[
(M_{p,q}^\lambda)' = \Bigl\{ f- \text{ measurable }: \sup_{\|g\|_{M_{p,q}^\lambda}=1} \int_{T^n} f(x) g(x) \, dx < \infty \Bigr\}.
\]

\begin{lemma} \label{lem2}\cite{ED&SD}  If $\lambda_0 \neq \lambda_1$, then the continuous embedding
\[
(M_{p,q}^\lambda)' \hookrightarrow \bigl( (M_{p,q_0}^{\lambda_0})', (M_{p,q_1}^{\lambda_1})' \bigr)_{\theta,q'}.
\]

\end{lemma}

\begin{lemma}\label{lem 8} If $1 \le p \le \infty$, $1 \le q \le \infty$ and $0 \le \alpha \le \frac{1}{p}$, then
\[
\|f\|_{(M_{p,q}^\alpha)'} \le \inf_\mathfrak{T} \Biggl( \sum_{m=0}^\infty \Bigl( 2^{m\alpha} \sum_{Q \in \mathfrak{T}_m = \mathfrak{T} \cap G_m} \|f\|_{L_{p'}(Q)} \Bigr)^{q'} \Biggr)^{1/q'} = \inf_\mathfrak{T} \|f\|_{LM_{p',q'}^{-\alpha}(\mathfrak{T})},
\]
where the infimum is taken over all dyadic partition $\mathfrak{T}$  of $\mathbb{T}^n$.  
\end{lemma}

\noindent
Proof of Lemma \ref{lem 8} .
By the definition of the dual space
\[
\|f\|_{(M_{p,q}^\alpha)'} = \sup_{\|g\|_{M_{p,q}^\alpha}=1} \int_{\mathbb{T}^n} f(x) g(x) \, dx<\infty.
\]
 Let $\mathfrak{T} = \{Q\}$ be an arbitrary dyadic partition of $\mathbb{T}^n$ . Therefore,
\begin{align*}
\left|\int_{\mathbb{T}^n} f(x) g(x) \, dx \right| &= \left| \sum_{m=0}^{\infty} \sum_{Q \in \mathfrak{T}_m} \int_Q f(x) g(x) \, dx \right| \\
&\le \sum_ {m=0}^{\infty} \sum_{Q \in \mathfrak{T}_m} \|f\|_{L_{p'}(Q)} \|g\|_{L_p(Q)} \\
&\le \sum_{m=0}^{\infty} \sum_{Q \in \mathfrak{T}_m} \|f\|_{L_{p'}(Q)} \sup_{Q \in G_m} \|g\|_{L_p(Q)} \\
&\le \left( \sum_{m=0}^{\infty} \Bigl( 2^{m\alpha} \sum_{Q \in \mathfrak{T}_m} \|f\|_{L_{p'}(Q)} \Bigr)^{q'} \right)^{1/q'} 
\left( \sum_{m=0}^{\infty} \Bigl( 2^{-m\alpha} \sup_{Q \in G_m} \|g\|_{L_p(Q)} \Bigr)^q \right)^{1/q} \\
&= \left( \sum_{m=0}^{\infty} \Bigl( 2^{m\alpha} \sum_{Q \in \mathfrak{T}_m} \|f\|_{L_{p'}(Q)} \Bigr)^{q'} \right)^{1/q'} \|g\|_{M_{p,q}^\alpha}.
\end{align*}

$\mathfrak {T}$ is an arbitrary dyadic partition, then
\[
\|f\|_{(M_{p,q}^\alpha)'} \le \inf_\mathfrak{T} \left( \sum_{m=0}^{\infty} \Bigl( 2^{m\alpha} \sum_{Q \in \mathfrak{T}_m} \|f\|_{L_{p'}(Q)} \Bigr)^{q'} \right)^{1/q'}.
\]
\endproof 
\begin{lemma}\label{lem 9}
Let $1 \le p \le \infty$, $0 \le \alpha, \lambda \le \frac{n}{p}$ and $\alpha + \lambda \ge \frac{n}{p}$, then
\begin{equation}\label{eq:lemma4}
\|f * g\|_{M_p^\lambda} \le c \, \|g\|_{(M_{p,\infty}^\alpha)'} \, \|f\|_{M_p^{\lambda+\alpha-\frac{n}{p}}},
\end{equation}
where $c$ is a constant depending only on $\lambda$, $n$ and $p$ .   
\end{lemma}

\noindent
Proof of Lemma \ref{lem 9} . 
Let  $g \in (M_{p,\infty}^\alpha)'$, $f \in M_p^{\lambda+\alpha-\frac{n}{p}}$.  A direct computation gives:
\begin{align*}
\|f*g\|_{M_p^\lambda} &= \sup_{0\le k \le\infty} 2^{-k\lambda} \sup_{\mathfrak{T} \in G_k} \|f*g\|_{L_p(\mathfrak{T})} \\
&\le \sup_{0\le k \le\infty} 2^{-k\lambda} \sup_{\mathfrak{T} \in G_k} \int_{\mathbb{T}^n} |g(y)| \, \|f(\cdot - y)\|_{L_p(\mathfrak{T})} \, dy \\
&\le \|g\|_{(M_{p,\infty}^\alpha)'} \sup_{0\le k \le\infty} 2^{-k\lambda} \sup_{\mathfrak{T} \in G_k} \sup_{0\le m \le\infty} \sup_{Q \in G_m} 2^{-m\alpha} \| \| f(x-y) \|_{L_p(\mathfrak{T})} \|_{L_p(Q)} \\
&= \|g\|_{(M_{p,\infty}^\alpha)'} \sup_{0\le k,m \le\infty} \sup_{\substack{\mathfrak{T} \in G_k \\ Q \in G_m}} 2^{-k\lambda} 2^{-m\alpha} \Biggl( \int_Q \int_\mathfrak{T} |f(x-y)|^p \, dx \, dy \Biggr)^{\frac{1}{p}}.
\end{align*}

 For $k \le m$, we have
$$
2^{-k\lambda} 2^{-m\alpha} \Biggl( \int_Q \int_\mathfrak{T} |f(x-y)|^p \, dx \, dy \Biggr)^{\frac{1}{p}} $$
$$\le 2^{-k\lambda} \Biggl( \int_\mathfrak{T} \bigl( \sup_{Q \in G_m} 2^{-m\alpha} \|f(x-y)\|_{L_p(Q)} \bigr)^p dx \Biggr)^{\frac{1}{p}}$$
$$ \le 2^{k(\frac{n}{p}-\lambda)}\cdot 2^{-m\alpha} \sup_{\mu Q=2^{mn}}\|f\|_{L_p(Q)} $$ $$\le2^{m(\frac{n}{p}-\alpha-\lambda)}\sup_{\mu Q=2^{mn}}\|f\|_{L_p(Q)}
\lesssim \|f\|_{M_p^{\lambda+\alpha-\frac{n}{p}}}.
$$
 On the other hand, for $m \le k$ we get

$$2^{-k\lambda} 2^{-m\alpha} \Biggl( \int_Q \int_\mathfrak{T} |f(x-y)|^p \, dx \, dy \Biggr)^{\frac{1}{p}} $$
$$\le 2^{-k\lambda} 2^{-m\alpha} \Biggl( \int_Q \sup_{\mathfrak{T} \in G_k} \int_\mathfrak{T} |f(x-y)|^p dx \Biggr)^{\frac{1}{p}}$$
$$ \le 2^{m(\frac{n}{p}-\alpha)}\cdot 2^{-k\lambda} \sup_{\mu \mathfrak{T}=2^{kn}}\|f\|_{L_p(\mathfrak{T})} $$
$$\le2^{k(\frac{n}{p}-\alpha)-k\lambda}\sup_{\mu \mathfrak{T}=2^{kn}}\|f\|_{L_p(\mathfrak{T})}
\lesssim \|f\|_{M_p^{\lambda+\alpha-\frac{n}{p}}}.$$

Thus, we have 
$$\|f*g\|_{M_p^{\lambda}}\le c \|g\|_{(M_{(p,\infty)}^\alpha)'}\|f\|_{M_p^{\lambda+\alpha-\frac{n}{p}}}.$$
With (\ref{eq:lemma4}), it completes the proof .\endproof
\begin{theorem}\label{th 2}
Let $\mathfrak{T}$ be some local partition of the space $\mathbb{T}^n$ .  
Let $0 < \max(q,1) \le p \le \infty, \quad 0 < \lambda < \frac{n}{q}, \quad 0 \le \gamma \le \frac{n}{p}, \quad 0 < \alpha = \gamma - \lambda + \frac{n}{q} < \frac{n}{p}$. If $f \in M_p^\gamma(\mathbb{T}^n)$, $g \in LM_{p',\infty}^{-\alpha}(\mathfrak{T})$, then  $f*g \in M_q^\lambda(\mathbb{T}^n)$ and the following inequality holds
\begin{equation}\label{eq:theorem2}
\|f*g\|_{M_q^\lambda(\mathbb{T}^n)} \le c \, \|f\|_{M_p^\gamma(\mathbb{T}^n)} \, \|g\|_{LM_{p',\infty}^{-\alpha}(\mathfrak{T})},
\end{equation}
where the constant $c$ depends only on the parameters $n,\lambda,q,\alpha,p$. 
\end{theorem}   

\noindent
Proof of Theorem \ref{th 2} . 1) Let us first prove the statement when $1 \le q = p < \infty$. Let $\alpha_0,\lambda_0,\alpha_1,\lambda_1$ be such that $0<\lambda_0<\lambda<\lambda_1<\frac{n}{p}$ and $$\alpha_0+\lambda_0 = \alpha_1 + \lambda_1 = \alpha+\lambda = \gamma + \frac{n}{p}.$$
According to Lemma \ref{lem 9} we have the following inequalities
\[
\|f*g\|_{M_{p,\infty}^{\lambda_i}} \le c_i \, \|f\|_{M_{p,\infty}^\gamma} \, \|g\|_{(M_{p,\infty}^{\alpha_i})'}, \quad i=0,1
\]
Then for the operator  $T g = f*g$ we obtain
\[
T : (M_{p,\infty}^{\alpha_i})' \to M_{p,\infty}^{\lambda_i}, 
\]
where
$$\|T\|_i \le c_i \, \|f\|_{M_{p,\infty}^\gamma}.$$
Aplying the real interpolation method, it yields \cite{ED&SD}

$$T : \bigl( (M_{p,\infty}^{\alpha_0})', (M_{p,\infty}^{\alpha_1})' \bigr)_{\theta,q} \to (M_{p,\infty}^{\lambda_0}, M_{p,\infty}^{\lambda_1})_{\theta,q},$$ 
here
$$\|T\| \le (c_0 \|f\|_{M_{p,\infty}^\gamma})^{1-\theta} (c_1 \|f\|_{M_{p,\infty}^\gamma})^\theta=c_0^{1-\theta}c_1^{\theta}\|f\|_{M_{p,\infty}^{\gamma}}.$$
Taking into account Theorem \ref{th7}, Lemma \ref{lem2} and Lemma \ref{lem 8}, we obtain
$$
\|f*g\|_{M_{p,q}^\lambda} \lesssim \|f*g\|_{{{(M_{p,\infty}^{\lambda_0}},{M_{p,\infty}^{\lambda_1}})}_{\theta,q}}$$
$$\le\|f\|_{M_{p,\infty}^\gamma} \, \|g\|_{((M_{p,\infty}^{\alpha_0})', (M_{p,\infty}^{\alpha_1})')_{\theta,q}}\lesssim \|f\|_{M_{p,\infty}^\gamma} \, \|g\|_{(M_{p,q'}^\alpha)'}$$ $$\le\|f\|_{M_{p,\infty}}\inf_{\mathfrak{T}}(\sum_{k\in \mathbb{Z}}(2^{k\alpha}\sum_{Q\in\mathfrak{T}_k}\|g\|_{L^{p'}(Q)})^q)^{\frac{1}{q}}$$ 
$$\le \|f\|_{M_{p,\infty}^\gamma} \, \|g\|_{LM_{p',q}^{-\alpha}(T)}.
$$
Hence, in particular for $q=\infty$ , we have
$$\|f*g\|_{M_{p,q}}^{\lambda}\le \|f\|_{M_{p,\infty}^{\gamma}}\|g\|_{LM_{p',q}^{-\alpha}(\mathfrak{T})}. $$
 Since $M_{p,\infty}^{\lambda}=M_{p}^{\lambda}$ , it gives  (\ref{eq:theorem2})  for $p=q$. 

2) Now let us consider the case  $0<\max(q,1)\le p\le \infty$ . According to $M_p^{\tilde \lambda} \hookrightarrow M_q^\lambda$ , where $\tilde \lambda = \lambda - \frac{1}{q} + \frac{1}{p}$, $q \le p$ (see Lemma \ref{lem 9} ) we obtain that
\[
\|f*g\|_{M_q^\lambda} \le \|f*g\|_{M_p^{\tilde \lambda}} \lesssim \|f\|_{M_p^\gamma} \, \|g\|_{LM_{p',q}^{-\alpha}(\mathfrak{T})}.
\]

Furthermore, by applying the inequality which was proved above we complete the proof.  \endproof 

\vspace{1em}

\textbf{Conclusion}

The study focuses on generalizing classical O’Neil-type inequalities for convolution operators in Morrey-type spaces defined on the torus $\mathbb{T}^n$. Relying on the fundamental principles of harmonic analysis, we derive new estimates for the mean behavior of functions. These results contribute to addressing important problems in the theory of differential equations and signal processing.
The work is aimed at extending the ideas of the classical Minkowski and Young inequalities, as well as O’Neil’s results in Lorentz spaces, thereby broadening their range of applications.
The main result of the study is the proof that, under certain conditions, the convolution $f * g$ belongs to a Morrey-type space with different parameters. This norm estimate is obtained using the properties of generalized Morrey and Köthe spaces, as well as interpolation methods. These results provide a deeper understanding of the properties of convolution operators.


\noindent

\selectlanguage{english}


\begin{thebibliography}{99}

\bibitem{Gl2014} Grafakos L., Classical Fourier Analysis. Third Edition, Graduate Texts in Math.,  Springer, New York, Number 249, 638 (2014).

\bibitem{Bennet&Sharley} Bennett, C. and Sharpley, R. (1988). Interpolation of Operators Pure and Applied Mathematics. Boston, MA, Academic Press, INC.,129, 469 (1988).

\bibitem{Conv 1963} O’Neil, R. Convolution operators and L(p,q) spaces. Duke Math. J., 30 (1), 129–142 (1963).

\bibitem{Tleu Sadyk19}  Tleukhanova N.T., Sadykova K.K., O’Neil-type inequalities for convolutions in anisotropic Lorentz spaces. Eurasian Mathematical Journal, 10 (3), 68-83(2019)

\bibitem{TleuSad2020}Tleukhanova N.T., Sadykova K.K., The convolution in anisotropic Triebel– Lizorkin spaces. Abai university Bulletin. Ser. Physics and Mathematical Sciences, 1 (69), 163-168 (2020)

\bibitem{NursTikh} E. Nursultanov, S. Tikhonov, Convolution inequalities in Lorentz spaces. Journal of Fourier Analysis and Applications, 17,  486–505 (2010)

\bibitem{ED&SD} Nursultanov E.D., Suragan D., On the convolution operator in Morrey spaces. Journal of Mathematical Analysis and Applications,  515(1),  Art. 126357, (2022)

\bibitem{Bur.Chigam} Burenkov V.I., Chigambayeva D.K., Nursultanov E.D., Marcinkiewicz-type interpolation theorem for Morrey-type spaces and its corollaries. Complex Variables and Elliptic Equations, 65(1), pp. 87–108(2020).

\end{thebibliography}
\end{document}